\RequirePackage[l2tabu, orthodox]{nag}
\RequirePackage{snapshot}

\documentclass[10pt,twocolumn]{article}

\makeatletter
\if@twocolumn
  \usepackage[dvips,letterpaper,top=0.5in, bottom=0.5in, left=0.75in, right=0.5in,includefoot,heightrounded]{geometry}
\else
  \usepackage[dvips,letterpaper,margin=1in,includefoot,heightrounded]{geometry}
\fi

\usepackage{srcltx}

\usepackage{amsmath}
\usepackage{amssymb,amsfonts}

\usepackage{abstract}

\usepackage{epsfig}
\usepackage[usenames,dvipsnames]{color}
\usepackage[usenames,dvipsnames]{xcolor}
\usepackage{subfigure}

\usepackage{booktabs}

\usepackage{setspace}
\usepackage{flushend}
\usepackage{multicol}

\usepackage{cite}
\usepackage{url}
\usepackage[normalem]{ulem}

\usepackage{enumerate}

\usepackage{hyperref}

\usepackage{lipsum}

\usepackage[sc,tiny]{titlesec}

\newtheorem{proposition}{Proposition}
\newtheorem{definition}{Definition}
\newtheorem{theorem}{Theorem}
\newtheorem{lemma}{Lemma}
\newtheorem{corollary}{Corollary}

\newcommand{\cas}{\operatorname{cas}}
\newcommand{\Sa}{\operatorname{Sa}}
\newcommand{\Ca}{\operatorname{Ca}}
\newcommand{\unidade}[2]{\ensuremath{#1\,\mathrm{#2}}}

\def\QED{\mbox{$\square$}}
\def\proof{\noindent{\it Proof:~}}
\def\endproof{\hspace*{\fill}~\QED\par\endtrivlist\unskip}



\title{%
A Short Survey on Arithmetic Transforms\\
and the Arithmetic Hartley Transform
}

\author{%
R.~J.~Cintra%
\thanks{%
R.~J.~Cintra
and
H.~M.~de Oliveira
were with
the Communications Research Group,
Departamento de Eletr\^onica e Sistemas,
Universidade Federal de Pernambuco, Recife, Brazil.
Currently
they are with
the Signal Processing Group,
Departamento de Estat\'{\i}stica
at the same university.
E-mail:
\protect\url{rjdsc@de.ufpe.br},
\protect\url{hmo@de.ufpe.br}
}
\and
H. M. de Oliveira%
}

\date{}

\newcommand{\myabstract}{%
Arithmetic complexity has a main role in the performance
of algorithms for spectrum evaluation.
Arithmetic transform theory offers a method for computing
trigonometrical transforms with minimal number of multiplications.
In this paper, the proposed algorithms for the arithmetic Fourier transform are
surveyed.
A new arithmetic transform
for computing the discrete Hartley transform is introduced:
the Arithmetic Hartley transform.
The interpolation process is shown to be
the key element
of the arithmetic transform theory.}

\newcommand{\mykeywords}{%
Arithmetic transforms,
discrete transforms,
Fourier series,
VLSI implementations
}

\begin{document}

\makeatletter
\if@twocolumn

\twocolumn[%
  \maketitle
  \begin{onecolabstract}
    \myabstract
  \end{onecolabstract}
  \begin{center}
    \small
    \textbf{Keywords}
    \\\medskip
    \mykeywords
  \end{center}
  \bigskip
]
\saythanks

\else

  \maketitle
  \begin{abstract}
    \myabstract
  \end{abstract}
  \begin{center}
    \small
    \textbf{Keywords}
    \\\medskip
    \mykeywords
  \end{center}
  \bigskip
  \onehalfspacing
\fi

\section{Introduction and Historical Background}

Despite the existence of fast algorithms for discrete transforms (e.g., fast Fourier transform, FFT),
it is well known that the number of multiplications
can significantly increase their computational (arithmetic) complexity.
Even today, the multiplication operation
consumes much more time than addition or subtraction.
Table~\ref{clockcount} brings the clock count of some mathematical operations
 as implemented for the
Intel Pentium{\small\texttrademark} processor.
Observe that multiplications and divisions can be by far more
time demanding than additions, for instance.
Sine and cosine function costs are also shown.

This fact stimulated the research on discrete transform algorithms
that minimize
the number of multiplications.
The Bhatnagar's algorithm~\cite{Bhatnagar:97},
which uses Ramanujan numbers to eliminate multiplications
(however, the choice of the transform blocklength is rather limited),
is an example.
Parallel to this, approximation approaches,
which perform a trade-off between accuracy and computational complexity,
have been proposed~\cite{Dee:01,Dimitrov:99,Rounded:02}.

\begin{table}
\centering
\caption{Clock count for some arithmetic instructions carried
on a Pentium{\small\texttrademark} processor~\cite{antonakos1996pentium}}

\label{clockcount}
\begin{tabular}{lc}
\toprule
Operation  & Clock count \\
\midrule
{\tt add} & 1--3 \\
{\tt sub} & 1--3 \\
{\tt fadd} & 1--7 \\
{\tt fsub} & 1--7 \\
{\tt mul} (unsigned) & 10--11 \\
{\tt mul} (signed) & 10--11 \\
{\tt div} (unsigned)& 17--41 \\
{\tt div} (signed) & 22--46 \\
{\tt fdiv} & 39 \\
{\tt sin}, {\tt cos} & 17--137 \\
\bottomrule
\end{tabular}
\end{table}

Arithmetic transforms emerged in this framework
as an algorithm for spectrum evaluation,
aiming the elimination of multiplications.
Thus, it may possess a lower computational complexity.
The theory of arithmetic transform is essentially based on M{\"o}bius function theorems~\cite{Burton},
offering only trivial multiplications, i.e., multiplications by $\{-1,0,1\}$.
Therefore, only addition operations (except for multiplications
by scale factors) are left to computation.
Beyond the computational attractiveness,
arithmetic transforms
turned out to be naturally suited for parallel processing
and VLSI design~\cite{Tufts:Sadasiv:88a,Reed:Shih:90}. %

\begin{figure}
\centering
\epsfig{file=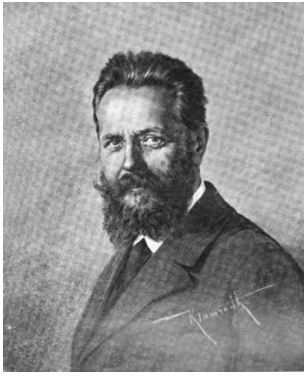,width=4cm}
\caption{Ernst Heinrich Bruns: pioneer of arithmetic transform theory.
Image in public domain.}
\label{figure-bruns}
\end{figure}

The very beginning ofthe research on arithmetic transforms dates back to 1903
when
the German mathematician
Ernst Heinrich Bruns\footnote{Bruns (1848-1919) earned his
doctorate in 1871 under the supervision of Weierstrass and Kummer.},
depicted in Figure~\ref{figure-bruns},
published the
\emph{Grundlinien des wissen\-schaft\-lichnen Rech\-nens}~\cite{Bruns:03},
the seminal work in this field.
In spite of that,
the technique remained unnoticed for a long time.
Forty-two years later, in Baltimore, \mbox{U.S.A.},
the Hungarian Aurel Freidrich Wintner%
\footnote{A curious fact: Wintner was born in April 8th, 1903, in Budapest,
the same year when Bruns had published his
\mbox{\emph{Grundlinien}}. Wintner died on January 15th, 1958, in Baltimore.},
privately published
a monograph entitled
{\em An Arithmetical Approach to Ordinary Fourier Series}.
This monograph presented an arithmetic method using M\"obius function to
calculate the Fourier series of even periodic functions.

After Wintner's monograph,
the theory entered again in ``hibernation'' state.
Not before 1988, Dr.~Donald W.~Tufts and Dr.~Angaraih G.~Sadasiv, independently,
had reinvented Wintner's arithmetical procedure,
reawaking the arithmetic transform.

In the quest to implement it,
two other researchers played an important role:
Dr.~Oved Shisha of the \mbox{U.R.I.}
Department of Mathematics and
Dr.~Charles Rader of Lincoln Laboratories.
They were aware of Wintner's monograph and helped Tufts
in many discussions.
In 1988
{\em The Arithmetic Fourier Transform} by Tufts and Sadasiv was published
in IEEE Acoustic, Speech, and Signal Processing (ASSP) Magazine~\cite{Tufts:Sadasiv:88a}.

Another breakthrough came in early 1990s
is due to
Emeritus Professor Dr.~Irving S.~Reed---%
the co-inventor of the
widely used
Reed-Muller~(1954) and
Reed-Solomon~(1964) codes.
Author of hundreds of publications, Dr.~Reed made
important contributions to the area of signal processing.
Specifically on arithmetic transforms,
in 1990
Reed,
Tufts and co-workers
provided two fundamental contributions~\cite{Reed:90,Reed:92}.
In~\cite{Reed:90},
a reformulated version of Tufts-Sadasiv approach,
the arithmetic Fourier transform (AFT) algorithm was
able to encompass a larger class of signals
and to compute
Fourier series of odd periodic functions as well as even periodic ones.
The publication of the 1992
{\em A VLSI Architecture for Simplified Arithmetic Fourier Transform Algorithm}
by Dr.~Reed and collaborators  in the IEEE Transactions on ASSP~\cite{Reed:92}
was another crucial slash on the subject.
Indeed,
that paper was previously presented at
the
{\em International Conference on Application Specific Array Processors}
held in Princeton.
However, the 1992 publication reached a vastly larger
public,
since it was published in a major journal.
The new method,
an enhancement of the last proposed
algorithm~\cite{Reed:90},
was re-designed to have a
more balanced and computationally efficient
performance.
As a matter of fact,
Reed~{\em et~al.} proved that the newly proposed algorithm
was
identical to Bruns' original method.

When the AFT was introduced, some concerns on the feasibility of the AFT
were pointed out~\cite{Tepedelenlioglu:89}.
The main issue dealt with the number of samples required by the algorithm.
However, later studies showed that the use of interpolation techniques
on a sub-sampled set (e.g., zero- and first-order interpolation)
could overcome these difficulties~\cite{Tufts:89}.

The conversion of the standard 1-D AFT
into 2-D versions was just a matter of time.
Many variants were proposed
following the same guidelines
of the 1-D case~\cite{Reed:Choi:89,
Choi:90,
Huisheng:96,
Atlas:97,
Ge:97,
Choi:89}.
Further research was carried out seeking different implementations of the AFT.
An alternative method~\cite{Lovine:93}
proposed a
``M\"obius-function-free AFT''.
Iterative~\cite{Tufts:93} and adaptative approaches~\cite{Li:90}
were also examined.
In spite of that, the most popular presentations of the AFT are still those found
in~\cite{Reed:90,Reed:92}.

Although the main and original motivation of the arithmetic algorithm
was the computation
of the Fourier Transform,
further generalizations were performed and the arithmetic approach was utilized to
calculate other transforms.
Dr.~Luc Knockaert of Department of Information Technology at
Ghent University, Belgium, amplified the Bruns procedure,
defining a generalized M\"obius transform~\cite{Knockaert:94,Knockaert:96}.
Moreover, four versions of the cosine transform was shaped
in the arithmetic transform formalism~\cite{Huisheng:Ping:96}.

Further generalization came in early 2000s with the definition
of the Arithmetic Hartley Transform~(AHT)~\cite{Cintra:Crete:02,Cintra:Florida:02}.
These works constituted an effort to make arithmetical procedure
applicable for the computation of trigonometrical transforms, other than Fourier transform.
In particular the AHT computes the discrete Hartley
transform\footnote{Ralph Vinton Lyon Hartley (1888-1970)
introduced his real integral transform in a 1942 paper published in the
{\em Proceedings of I.R.E.}
The Hartley transform relates a pair of signals $f(t)\longleftrightarrow F(\nu)$ by
\begin{eqnarray*}
F(\nu) = \frac{1}{\sqrt{2\pi}}\int_{-\infty}^{\infty}f(t)(\cos(\nu t) + \sin(\nu t))\mathrm{d}t, \\
f(t) = \frac{1}{\sqrt{2\pi}}\int_{-\infty}^{\infty}F(\nu)(\cos(\nu t) + \sin(\nu t))\mathrm{d}\nu.
\end{eqnarray*}
}:
the real, symmetric, Fourier-like discrete transform defined in 1983
by Emeritus Professor Ronald Newbold Bracewell
in \textsl{The Discrete Hartley Transform}, an article published in
the Journal of Optical Society of America.

In 1988 and then the technological state-of-art was
dramatically different from that Bruns and Wintner found.
Computational facilities and digital signal processing
integrated circuits made possible AFT to leave theoretical
constructs and reach practical implementations.
Since its inception in engineering,
the AFT was recognized as tool to be
implemented with VLSI techniques.
Tufts himself had observed that AFT could be naturally implemented in VLSI architectures~\cite{Tufts:Sadasiv:88a}.
Implementations were proposed in
\cite{Fischer:90,
Wigley:Jullien:90,
Wigley:Jullien:92,
Park:91,
Park:Prassana:91,
Kelley:93,
Park:Prassana:93,
Reed:Shih:90,
Reed:92,
Atlas:97,
Fischer:89,
Park:Prassana:tech:90,
Park:Prassana:tech:91,
DiLecce:95}.
Initial applications of the AFT took place in several areas:
pattern matching techniques~\cite{Zaid:92},
measurement and instrumentation~\cite{Andria:96,Andria:DiLecce:96},
auxiliary tool for computation of $z$-transform~\cite{Hsu:Reed:94,Hsu:94},
and imaging~\cite{Tufts:Fan:89}.

This paper is organized in two parts.
In Section~\ref{aft}, the mathematical evolution of the Arithmetic Fourier Transform is outlined.
In Section~\ref{aht}, a summary of the major results on the Arithmetic Hartley Transform is shown.
Interpolation issues are addressed and many points of the AFT are clarified,
particularly the zero-order approximation.

\section{The Arithmetic Fourier Transform}
\label{aft}

In this section,
the three major
developments of the arithmetic Fourier transform technique are presented.
With emphasis on the theoretical groundwork,
the AFT algorithms devised by Tufts, Sadasiv, Reed {\em et alli}
are briefly outlined.
In this work,
$k_1|k_2$ denotes that~$k_1$ is a divisor of~$k_2$;
$\lfloor \cdot \rfloor$ is the floor function and~$[\cdot]$ is the nearest integer function.

\begin{lemma}\label{soma:cos:sin}
Let $k$, $k'$ and $m$ be integers.
\begin{equation}
\sum_{m=0}^{k-1}\cos\left({2\pi m\frac{k'}{k}}\right)=
\begin{cases}
k& \text{if $k|k'$}, \\
0& \text{otherwise} \cr
\end{cases}
\end{equation}
and
\begin{equation}
\sum_{m=0}^{k-1}\sin\left({2\pi m\frac{k'}{k}}\right)=0.
\end{equation}
\end{lemma}
\proof
Consider the expression $\sum_{m=0}^{k-1}\left({\rm e}^{2\pi j\frac{k'}{k}}\right)^m$.
When $k|k'$, yields
\[
\sum_{m=0}^{k-1}\left({\rm e}^{2\pi j\frac{k'}{k}}\right)^m=\sum_{m=0}^{k-1}1=k.
\]
Otherwise,
we have:
$$
\sum_{m=0}^{k-1}\left({\rm e}^{2\pi j\frac{k'}{k}}\right)^m=
\frac{1-{\rm e}^{j2\pi k'}}{ 1-{\rm e}^{j2\pi \frac{k'}{k}}}=0.
$$
Therefore,
it follows that:
$$
\sum_{m=0}^{k-1}{\rm e}^{2\pi jm\frac{k'}{k}}=
\begin{cases}
k& \text{if $k|k'$}, \\
0& \text{otherwise}.
\end{cases}
$$
By taking the real and imaginary parts,
we conclude the proof.
\endproof

\begin{definition}[M{\"o}bius $\mu$-function]
For a positive integer~$n$,
\begin{equation}
\mu(n)\triangleq
\begin{cases}
1& \text{if $n=1$}, \\
{(-1)}^r& \text{if $n=\prod_{i=1}^rp_i$, $p_i$ distinct primes}, \\
0& \text{if $p^2|n$ for some prime $p$}.
\end{cases}
\end{equation}
\end{definition}

A relevant lemma related to $\pmb{\mu}$-function is stated below.
\begin{lemma}\label{soma:mu}
\begin{equation}
\sum_{d|n} \mu(d)=
\begin{cases}
1& \text{if $n=1$}, \\
0& \text{if $n>1$}.
\end{cases}
\end{equation}
\end{lemma}

\begin{theorem}[M{\"o}bius Inversion Formula for Finite Series]
\label{reed:2}
Let $n$ be a positive integer and  $f_n$ a non-null sequence
for $1\leq n\leq N$ and null for $n>N$.
If
\begin{equation}
g_n=
\sum_{k=1}^{\lfloor N/n \rfloor}f_{kn},
\end{equation}
then
\begin{equation}
f_{n}=
\sum_{m=1}^{\lfloor N/n \rfloor}\mu(m)g_{mn}.
\end{equation}
\endproof
\end{theorem}
This is the finite version of the
M{\"o}bius inversion formula~\cite{Burton}.
A proof can be found in~\cite{Reed:90}.

\subsection{Tufts-Sadasiv Approach}

Consider a real even periodic function expressed by its Fourier series, as seen below:
\begin{equation}\label{sadasiv:5}
v(t)=\sum_{k=1}^{\infty}v_k(t).
\end{equation}
The components $v_k(t)$ represent the harmonics of $v(t)$, given by:
\begin{equation}\label{sadasiv:6}
v_k(t)=a_k\cdot\cos(2\pi kt),
\end{equation}
where $a_k$
is the amplitude
of the $k$th harmonic.

It was assumed, without loss of generality, that $v(t)$
had unitary period and null mean ($a_0=0$).
Furthermore, consider the $N$ first harmonics as the only significant ones,
in such a way that $v_k(t)=0$, for $k>N$ (bandlimited approximation).
Thus the summation of Equation~\ref{sadasiv:5} might be constrained to $N$~terms.

\begin{definition}\label{def:avg:sadasiv}
The $n$th average is defined by
\begin{equation}\label{sadasiv:8a}
S_n(t)\triangleq\frac{1}{n}\sum_{m=0}^{n-1}v\left(t-\frac{m}{n}\right),
\end{equation}
for $n=1,2,\ldots,N$.
$S_n(t)$ is null for $n>N$.
\endproof
\end{definition}

After an application of Equations~\ref{sadasiv:5}
and~\ref{sadasiv:6} into~\ref{sadasiv:8a},
it yielded:
\begin{align}
S_n(t)=&\frac{1}{n}\sum_{m=0}^{n-1}v\left(t-\frac{m}{n}\right) \nonumber \\
=&\frac{1}{n}\sum_{m=0}^{n-1}\sum_{k=1}^{\infty}a_{k}\cos\left(2\pi kt-2\pi k\frac{m}{n}\right) \nonumber \\
=&\frac{1}{n}\sum_{k=1}^{\infty}a_{k}\sum_{m=0}^{n-1}\left(\cos(2\pi kt)\cos\left(2\pi k\frac{m}{n}\right)\right. - \nonumber \\
& \left. \sin(2\pi kt)\sin\left(2\pi k\frac{m}{n}\right)\right) \nonumber \\
=&\frac{1}{n}\sum_{k=1}^{\infty}a_{k}\cos(2\pi kt)\cdot \left.
\begin{cases}
n&\text{if $n|k$}, \\ 0& \text{otherwise}
\end{cases}
\right\} \nonumber \\
=&\sum_{n|k}^{\infty}v_{k}(t)=\sum_{m=1}^{\infty}v_{mn}(t),
\quad n=1,\ldots,N.\label{sadasiv:9}
\end{align}
Proceeding that way, the $n$th average can be written in terms of
the harmonics of $v(t)$, instead of its samples (Definition~\ref{def:avg:sadasiv}).
Since we assumed
$v_n(t)=0, n>N$,
only the first $\lfloor N/n\rfloor$
terms of Equation~\ref{sadasiv:9} might possibly be nonnull.

As a consequence the task was to invert Equation~\ref{sadasiv:9}.
Doing so, the harmonics can be expressed in terms of the averages, $S_n(t)$,
which are derived from the samples of the signal $v(t)$.
The inversion is accomplished by invoking the M\"obius inversion formula.

\begin{theorem}
The harmonics of $v(t)$ can be obtained by:
\begin{equation}\label{sadasiv:11}
v_{k}(t)=\sum_{m=1}^{\infty}\mu(m)S_{mk}(t),\quad \forall k=1,\ldots, N.
\end{equation}
\end{theorem}
\proof
Some manipulation is needed.
Substituting Equation~\ref{sadasiv:9} in Equation~\ref{sadasiv:11},
it yields
\begin{equation}
\sum_{m=1}^{\infty}\mu(m)S_{mk}(t)=
\sum_{m=1}^{\infty}\mu(m)\sum_{n=1}^{\infty}v_{kmn}(t).
\end{equation}
Now it is the tricky part of the proof.
\begin{equation}
\begin{split}
\sum_{m=1}^{\infty}\mu(m)\sum_{n=1}^{\infty}v_{kmn}(t)
&=\sum_{m=1}^{\infty}\sum_{n=1}^{\infty}\mu(m)v_{kmn}(t) \\
&=\sum_{j=1}^{\infty}v_j(t)\left( \sum_{m| \frac{j}{k}}^{\infty}\mu(m) \right).%
\end{split}
\end{equation}
According to Lemma~\ref{soma:mu},
the inner summation can only be null if $j/k=1$.
In other words, the term $v_k(t)$ is the only survivor of the outer summation and
the proof is completed.\endproof

The following aspects of the Reed-Tufts algorithm can be
highlighted~\cite{Tufts:Sadasiv:88a}:
\begin{itemize}
\item
This initial version of the AFT had a strong constraint:
it can only handle even signals;

\item
All computations are performed using only additions (except for few
multiplications due to scaling);

\item
The algorithm architecture is suitable for parallel processing, since
each average is computed independently from the others;

\item
The arithmetic transform theory is based on Fourier series, instead of the discrete Fourier transform.
\end{itemize}

\subsection{Reed-Tufts Approach}%

Presented by Reed~\emph{et al.} in 1990~\cite{Reed:90},
this algorithm is a generalization of Tuft-Sadasiv method.
The main constraint of the latter procedure
(handling only with even signals) was removed,
opening path for the computation of all Fourier series coefficient of
periodic functions.

Let $v(t)$ be a real $T$-periodic function, whose $N$-term finite Fourier series
is given by
\begin{equation}
v(t)=a_0 +
\sum_{n=1}^{N}a_n\cos\left( \frac{2\pi nt}{T}\right)+
\sum_{n=1}^{N}b_n\sin\left( \frac{2\pi nt}{T}\right),
\end{equation}
where $a_0$ is the mean value of $v(t)$.
The even and odd coefficients of the Fourier series are
$a_n$ and $b_n$, respectively.

Let $\bar{v}(t)$ denote the signal~$v(t)$ removed of its mean value~$a_0$.
Consequently,
\begin{equation}
\begin{split}
\bar{v}(t)&=v(t)-a_0 \\
&=\sum_{n=1}^{N}a_n\cos\left( \frac{2\pi nt}{T}\right)+
\sum_{n=1}^{N}b_n\sin\left( \frac{2\pi nt}{T}\right).
\end{split}
\end{equation}
A delay (shift) of $\alpha T$ in $\bar{v}(t)$ furnishes to the following
expression:
\begin{equation}
\label{reed:6a}
\begin{split}
\bar{v}(t+\alpha T)
=&\sum_{n=1}^{N}a_n \cos\left(2\pi n(\frac{t}{T}+\alpha)\right)+ \\
 &\sum_{n=1}^{N}b_n \sin\left(2\pi n(\frac{t}{T}+\alpha)\right)  \\
=&\sum_{n=1}^{N}c_n(\alpha) \cos\left(2\pi n\frac{t}{T}\right) + \\
 &\sum_{n=1}^{N}d_n(\alpha) \sin\left(2\pi n\frac{t}{T}\right),
\end{split}
\end{equation}
where $-1<\alpha<1$ and
\begin{align}
c_n(\alpha)&=a_n \cos(2\pi n\alpha)+b_n\sin(2\pi n\alpha), \label{reed:6b}\\
d_n(\alpha)&=-a_n \sin(2\pi n\alpha)+b_n\cos(2\pi n\alpha).
\end{align}

In the sequel,
the computation of the Fourier coefficients~$a_n$ and $b_n$
based on~$c_n(\alpha)$ is outlined.
Meanwhile, the formula for the $n$th average
(Tufts-Sadasiv)
is updated by the next definition.
\begin{definition}\label{def:soma_parcial}
The $n$th average is given by
\begin{equation}\label{reed:8}
S_n(\alpha)\triangleq\frac{1}{n}\sum_{m=0}^{n-1}\bar{v}\left(\frac{m}{n}T+\alpha T\right),
\end{equation}
where $-1< \alpha <1$.
\end{definition}

Now the quantities $c_n(\alpha)$ can be
related to the averages, according to the following theorem.
\begin{theorem}\label{theor:c_n}
The coefficients $c_n(\alpha)$ are computed via
M\"obius inversion formula for finite series
and are expressed by
\begin{equation}
c_n(\alpha)=
\sum_{l=1}^{\lfloor N/n \rfloor}\mu(l)S_{ln}(\alpha).
\end{equation}
\end{theorem}
\proof
Substituting the result of Equation~\ref{reed:6a} into Equation~\ref{reed:8}
furnishes the following expression:
\begin{equation}
\begin{split}
S_n(\alpha)=&
\sum_{k=1}^{N}c_{k}(\alpha)\frac{1}{n}\sum_{m=0}^{n-1}\cos \left(\frac{2\pi km}{n}\right)
+ \\
& \sum_{k=1}^{N}d_{k}(\alpha)\frac{1}{n}\sum_{m=0}^{n-1}\sin \left(\frac{2\pi km}{n}\right).
\end{split}
\end{equation}
A direct application of Lemma~\ref{soma:cos:sin} yields
\begin{equation}
S_n(\alpha)=\sum_{l=1}^{\lfloor N/n \rfloor} c_{ln}(\alpha).
\end{equation}
Invoking the
M\"obius inversion formula for finite series,
the theorem is proved.
\endproof

\noindent
Finally, the main result can be derived.
\begin{theorem}[Reed-Tufts]
\label{theor:reed-tufts}
The Fourier series coefficients $a_n$ and $b_n$ are computed by
\begin{align}
a_n&=c_n(0), \label{reed:13}\\
b_n&=(-1)^mc_n\left(\frac{1}{2^{k+2}}\right) \quad n=1,\ldots,N,
\end{align}
where $k$ and $m$ are determined by the factorization $n=2^k(2m+1)$.
\end{theorem}
\proof
For $\alpha=0$, using Equation~\ref{reed:6b},
it is straightforward to show that $a_n=c_n(0)$.
For $\alpha=\frac{1}{2^{k+2}}$ and $n=2^k(2m+1)$, there are two sub-cases:
$m$ even or odd.
\begin{itemize}
\item For $m=2q$, $n=2^k(4q+1)$.
Therefore,
\begin{equation}
2\pi n\alpha = 2\pi \frac{2^k(4q+1)}{2^{k+2}}= 2\pi q +\frac{\pi}{2}.
\end{equation}
Consequently, substituting this quantity in Equation~\ref{reed:6b}, yields
\begin{equation}
\begin{split}
c_n\left(\frac{1}{2^{k+2}}\right)&=a_n \cos\left(2\pi q +\frac{\pi}{2}\right)+b_n\sin\left(2\pi q +\frac{\pi}{2}\right) \\
&=b_n.
\end{split}
\end{equation}
\item For $m=2q+1$, $n=2^k(4q+3)$.
It follows that
\begin{equation}
2\pi n\alpha = 2\pi \frac{2^k(4q+3)}{2^{k+2}}= 2\pi q +\frac{3\pi}{2}.
\end{equation}
Again invoking the Equation~\ref{reed:6b}, the following expression is derived.
\begin{equation}
\begin{split}
c_n\left(\frac{1}{2^{k+2}}\right)=&a_n \cos\left(2\pi q +\frac{3\pi}{2}\right)+\\
& b_n\sin\left(2\pi q +\frac{3\pi}{2}\right) \\
=&-b_n.
\end{split}
\end{equation}
\end{itemize}
Joining these two sub-cases, it is easy to verify that
\begin{equation}
b_n=(-1)^mc_n\left(\frac{1}{2^{k+2}}\right).
\end{equation}\endproof

The number of real multiplications and additions of this algorithm is
given by~\cite{Reed:90}
\begin{align}
M_R(N)&=\frac{3}{2}N, \\
\intertext{and}
A_R(N)&=\frac{3}{8}N^2,
\end{align}
respectively,
where $N$ is the blocklength of the transform.

\subsection{Reed-Shih (Simplified AFT)}

Introduced by Reed~\emph{et al.}~\cite{Reed:Shih:90},
this algorithm is an evolution of that one developed
by Reed and Tufts.
Surprisingly, in this new method, the averages are re-defined
in accordance to
the theory created by H.~Bruns~\cite{Bruns:03} in 1903.

\begin{definition}[Bruns Alternating Average]
The $2n$th Bruns alternating average, $B_{2n}(\alpha)$,
is defined by
\begin{equation}
B_{2n}(\alpha)\triangleq\frac{1}{2n}
\sum_{m=0}^{2n-1}(-1)^m\cdot v\left({\frac{m}{2n}T+\alpha T}\right).
\end{equation}
\end{definition}

Invoking the definition of $c_n$, applying
Theorem~\ref{theor:c_n} and Definition~\ref{def:soma_parcial},
the following theorem was derived.

\begin{theorem}\label{coef:bruns}
The coefficients $c_n(\alpha)$ are given by the
M{\"o}bius inversion formula for finite series as
\begin{equation}
c_n(\alpha)=\sum_{l=1,3,\ldots}^{\lfloor \frac{N}{n}\rfloor}
\mu(l)\cdot B_{2nl}(\alpha).
\end{equation}
\end{theorem}
\proof
See~\cite{Reed:92}.
\endproof

Since a relation between the signal samples and
the Bruns alternating averages was obtained, as well as
an expression connecting
the Bruns alternating averages to the $c_n$ coefficients,
was available,
few points are missing to compute the Fourier series coefficients.
Actually,
an expression relating
the Fourier series coefficients ($a_n$ and $b_n$)
to
the coefficients $c_n$
is sought.
Examining Equation~\ref{reed:6b}, two conditions
are distinguishable:
\begin{itemize}
\item $a_n=c_n(0)$;
\item $b_n=c_n\left( \frac{1}{4n}\right)$.
\end{itemize}
Calling Theorem~\ref{coef:bruns}, the next result was obtained.

\begin{theorem}[Reed-Shih]\label{coef:fourier-bruns}
The Fourier series coefficients  $a_n$ and $b_n$ are computed by
\begin{align}
a_0&=\frac{1}{T}\int_{0}^{T}v(t)\mathrm{d}t, \\
a_n&=\sum_{l=1,3,5,\ldots}^{\lfloor \frac{N}{n}\rfloor}\mu(l)B_{2nl}(0), \\
b_n&=\sum_{l=1,3,5,\ldots}^{\lfloor \frac{N}{n}\rfloor}\mu(l)(-1)^\frac{l-1}{2}B_{2nl}\left(\frac{1}{4nl}\right),
\end{align}
for $n=1,\ldots,N$.
\end{theorem}
\proof
The proof is similar to the proof of Theorem~\ref{theor:reed-tufts}.
\endproof

For a blocklength~$N$, the multiplicative and additive complexities
are given by
\begin{align}
M_R(N)&=N, \\
\intertext{and}
A_R(N)&=\frac{1}{2}N^2,
\end{align}
respectively.

The AFT algorithm proposed by Reed-Shih presents
some improvements over previous algorithms:
\begin{itemize}
\item
The computation of both $a_n$ and $b_n$
has
roughly the same computational effort.
The algorithm is more ``balanced'' than Reed-Tufts algorithm;

\item
The algorithm is naturally suited to a parallel processing
implementation;

\item
It is computationally less complex that Reed-Tufts algorithm.
\end{itemize}

\subsection{An Example}

In this subsection,
we draw some comments
in connection to
an example originally proposed by Reed~\emph{et al.}~\cite{Reed:Shih:90}.
Let $v(t)$ be a signal with period $T=\unidade{1}{s}$.
Consider the computation of the Fourier series coefficients
up to the 5th harmonic.

According to the Reed-Shih algorithm, the coefficients $a_n$ and
$b_n$ of the Fourier series of $v(t)$ are expressed by
\begin{equation}
\begin{bmatrix}
a_1\\ a_2\\ a_3\\ a_4\\ a_5
\end{bmatrix}=
\begin{bmatrix}
1& 0& -1& 0& -1 \\
0& 1 & 0& 0& 0 \\
0& 0& 1& 0& 0 \\
0& 0& 0& 1& 0 \\
0& 0& 0& 0& 1
\end{bmatrix}
\left[\begin{array}{l}
B_2(0) \\ B_4(0) \\ B_6(0) \\ B_8(0) \\ B_{10}(0)
\end{array}\right]
\end{equation}
and
\begin{equation}
\begin{bmatrix}
b_1\\ b_2\\ b_3\\ b_4\\ b_5
\end{bmatrix}=
\begin{bmatrix}
1& 0& 1& 0& -1 \\
0& 1 & 0& 0& 0 \\
0& 0& 1& 0& 0 \\
0& 0& 0& 1& 0 \\
0& 0& 0& 0& 1
\end{bmatrix}
\left[\begin{array}{l}
B_2(\frac{1}{4}) \\
B_4(\frac{1}{8}) \\
B_6(\frac{1}{12}) \\
B_8(\frac{1}{16}) \\
B_{10}(\frac{1}{20})
\end{array}\right].
\end{equation}
Comparing these formulations with the ones of Reed-Tufts algorithm,
one may
note the balance in the computation of $a_n$ and $b_n$.
Both coefficients are obtained through similar matrices.
Table~\ref{tabela:amostras} relates
the Bruns alternating averages $B_n(\alpha)$
to their required samples.
Notice that
at least 40~non-uniform time samples of $v(t)$ are necessary
to exactly compute the Bruns alternating averages.

\begin{table}[h]
\begin{center}
\caption
{Necessary samples for the Bruns alternating averages}
\label{tabela:amostras}
\begin{tabular}{cl}

\toprule
Bruns averages         & Sample time (s)\\
\midrule
$B_2(0)$ & $0, \frac{1}{2}$ \\
\noalign{\medskip}
$B_4(0)$ & $0, \frac{1}{4}, \frac{1}{2}, \frac{3}{4}$ \\
\noalign{\medskip}
$B_6(0)$ & $0, \frac{1}{6}, \frac{1}{3}, \frac{1}{2}, \frac{2}{3}, \frac{5}{6}$ \\
\noalign{\medskip}
$B_8(0)$ & $0, \frac{1}{8}, \frac{1}{4}, \frac{3}{8}, \frac{1}{2}, \frac{5}{8}, \frac{3}{4}, \frac{7}{8}$ \\
\noalign{\medskip}
$B_{10}(0)$ & $0, \frac{1}{10}, \frac{1}{5}, \frac{3}{10}, \frac{2}{5}, \frac{1}{2}, \frac{3}{5}, \frac{7}{10}, \frac{4}{5}, \frac{9}{10}$ \\
\noalign{\medskip}
$B_2(\frac{1}{4})$     & $\frac{1}{4},\frac{3}{4}$\\
\noalign{\medskip}
$B_4(\frac{1}{8})$     & $\frac{1}{8},\frac{3}{8},\frac{5}{8},\frac{7}{8}$\\
\noalign{\medskip}
$B_6(\frac{1}{12})$    & $\frac{1}{12},\frac{1}{4},\frac{5}{12},\frac{7}{12},\frac{3}{4},\frac{11}{12}$\\
\noalign{\medskip}
$B_8(\frac{1}{16})$    & $\frac{1}{16},\frac{3}{16},\frac{5}{16},\frac{7}{16},\frac{9}{16},\frac{11}{16},\frac{13}{16},\frac{15}{16}$\\
\noalign{\medskip}
$B_{10}(\frac{1}{20})$ & $\frac{1}{20},\frac{3}{20},\frac{1}{4},\frac{7}{20},\frac{9}{20},\frac{11}{20},\frac{13}{20},\frac{3}{4},\frac{17}{20},\frac{19}{20}$ \\
\bottomrule
\end{tabular}
\end{center}
\end{table}

At this point, some observations are relevant:
\begin{itemize}
\item
This algorithm is not naturally suited for
uniform sampling.
\item
A uniform sampler utilized
to obtain all the required samples
would need a very high sampling rate.
In the example illustrated here,
a \unidade{120}{Hz} clock should be
required to sample the necessary points for the computation of the Fourier series
of a \unidade{1}{Hz} bandlimited signal.
\end{itemize}

Certainly these observations appear to be disturbing and seems to jeopardize the feasibility of the whole procedure.
However, it is important to stress that this procedure
furnishes the {\em exact} computation of the Fourier series
coefficients.

An empirical solution to circumvent this problem is to interpolate.
An interpolation based on uniformly sampled points can be used to
estimate the sample values required by AFT.
Of course, this procedure inherently introduces computation errors.

For example, assuming that the 1~Hz signal $v(t)$
is sampled by a clock with period $T_0=\unidade{\frac{1}{10}}{s}$.
Hence, the following sample points are available:
\begin{align*}
v\left( 0\right),
v\left( \frac{1}{10}\right),
v\left( \frac{2}{10}\right),
v\left( \frac{3}{10}\right),
v\left( \frac{4}{10}\right), \\
v\left( \frac{5}{10}\right),
v\left( \frac{6}{10}\right),
v\left( \frac{7}{10}\right),
v\left( \frac{8}{10}\right),
v\left( \frac{9}{10}\right).
\end{align*}

Table~\ref{tabela:amostras} shows, for example, that the computation of
$B_4(0)$ requires  --- among other samples --- $v\left( \frac{1}{4}\right)$,
which is clearly not available.
To overcome this difficulty,
a rounding operation can be introduced.
Thus, the sample $v\left( \frac{3}{10}\right)$
can be used
whenever the algorithm called $v\left( \frac{1}{4}\right)$ ($\left[10 \frac{1}{4} \right]/10 = 3/10$).
This rounding operation is also known as zero-order interpolation.

The accuracy of the AFT algorithm is deeply associated
with the sampling period $T_0$.
If more precision is required, then one should expect to
increase sampling rate, resulting in the introduction of smaller errors
due to the interpolation scheme.
Higher order of interpolation (e.g. first-order interpolation) can also be used to
obtain more accurate estimations of the Fourier
series coefficients.
The following trade-off is quite clear
accuracy versus order of interpolation.

However, for signals sampled at Nyquist rate (or close to),
zero-order interpolation already leads to good results~\cite{Reed:90}.
A detailed error analysis of interpolation schemes
can be found in~\cite{Bartels:89,Reed:90,Huisheng:Ping:96,Hsu:Reed:94}.
Further comments can be found in~\cite{Reed:Shih:90}.

\section{A New Arithmetic Transform}%
\label{aht}

Besides its numerical appropriateness~\cite{Bracewell},
the discrete Hartley transform (DHT)
has proved along the years to be an important tool
with several applications, such as
biomedical image compression,
OFDM/CDMA systems, and
ADSL transceivers.
Searching the literature, no mention about
a possible ``Arithmetic Hartley Transform'' to compute
the DHT was found.

In this section, a condensation of the
main results of the Arithmetic Hartley Transform
is outlined.
The method used to define the AHT turned out to
furnish a new insight into the arithmetic transform.
In particular, the role of interpolation is clarified.
Additionally, it is mathematically shown
that interpolation is a pivotal issue
in arithmetic transforms.
Indeed it determines the transform.

A new approach to arithmetic transform is adopted.
Instead of considering uniformly sampled points extracted from
a continuous signal~$v(t)$,
the AHT is based on the purely discrete signal.
Thus, the starting point of the development is the discrete transform definition,
not the series expansion, as it was done in the development of the AFT algorithm.
This approach is philosophically appealing, since
in a final analysis a discrete transform relates two set of points,
not continuous functions.

Let $\mathbf{v}$ be an $N$-dimensional vector with real elements.
The DHT establishes a pair denoted by
$
\mathbf{v} = [v_0, v_1, \ldots, v_{N-1}]^\top
\leftrightarrow
\mathbf{V} = [V_0,V_1, \ldots, V_{N-1}]^\top
$,
where the elements of the transformed vector~$\mathbf{V}$
(i.e., Hartley spectrum)
are defined by~\cite{Bracewell}
\begin{equation}\label{DHT}
V_k
\triangleq
\frac{1}{N}
\sum_{i=0}^{N-1}
v_i\cdot
\cas \left(  \frac{2\pi k i}{N} \right),
\quad
k=0, 1, \ldots,N-1,
\end{equation}
where $\cas x \triangleq \cos x +\sin x$ is Hartley's
``cosine and sine'' kernel. %
The inverse discrete Hartley transform is then~\cite{Bracewell}
\begin{equation}\label{IDHT}
v_i
=
\sum_{k=0}^{N-1}
V_k\cdot
\cas \left(  \frac{2\pi k i}{N} \right),
\quad
i=0, 1, \ldots,N-1.
\end{equation}

\begin{lemma}[Fundamental Property]
\label{lema_fundamental}
The function $\cas(\cdot)$ satisfies
\begin{equation}\label{fundamental}
\sum_{m=0}^{k-1}\cas\left(2\pi m\frac{k'}{k}\right)=
\begin{cases}
k& \text{if $k|k'$}, \\
0& \text{otherwise}.
\end{cases}
\end{equation}
\proof
It follows directly from Lemma~\ref{soma:cos:sin}.
\endproof
\end{lemma}

Similarly to the AFT theory,
it is necessary to define averages $S_k$, calculated from the
time-domain elements. The averages are computed by
\begin{equation}
\label{soma:das:amostras}
S_k\triangleq\frac{1}{k}\sum_{m=0}^{k-1}v_{m\frac{N}{k}},
\qquad k=1,  \ldots, N-1.
\end{equation}

This definition
requires fractional index sampling.
Analogously to the AFT methods, this fact seems
to make further considerations impracticable,
since only integer index samples are available.
This subtle question is to be treated in the sequel.
For now,
we assume that the fractional index sample
are \emph{somehow} available.

An application of the inverse Hartley transform (Equation~\ref{IDHT})
in Equation~\ref{soma:das:amostras} offered:
\begin{equation}\label{inversa}
S_k
=
\frac{1}{k}\sum_{m=0}^{k-1}\sum_{k'=0}^{N-1}
V_{k'}\cas\left(\frac{2\pi k'\frac{mN}{k}}{N}\right).
\end{equation}
Rearranging the summation order, simplifying,
 and calling Lemma~\ref{lema_fundamental}, it yielded:
\begin{equation}\label{somaV}
\begin{split}
S_k&=\frac{1}{k}
\sum_{k'=0}^{N-1}V_{k'}
\sum_{m=0}^{k-1}
\cas
\biggr(
2\pi m\frac{k'}{k}
\biggr) \\
&=\sum_{s=0}^{\lfloor(N-1)/k \rfloor}V_{sk}.
\end{split}
\end{equation}

For simplicity and without loss of generality, consider
a signal $\mathbf{v}$ with zero mean value, i.e.,
$\frac{1}{N}\sum_{i=0}^{N-1}v_i=0$.
Clearly, this consideration has no influence on the values of
$V_k,\ k\neq0$.
An application of the
modified M{\"o}bius inversion formula for finite series~\cite{Reed:90}
is sufficient to
obtain the final theorem to derive
the Arithmetic Hartley Transform.
According to Theorem~\ref{reed:2},
the result below follows.
\begin{theorem}[Reed {\em et alli}]\label{inv}
If
\begin{equation}\label{equacao:48}
S_k=\sum_{s=1}^{\lfloor (N-1)/k\rfloor}V_{sk}, \quad 1\leq k\leq {N-1},
\end{equation}
then
\begin{equation}\label{somaS}
V_k=\sum_{l=1}^{\lfloor(N-1)/k \rfloor}\mu(l)S_{kl},
\end{equation}
where $\mu(\cdot)$
is  M{\"o}bius function.\endproof
\end{theorem}

\begin{figure}[t]
\centering
\input{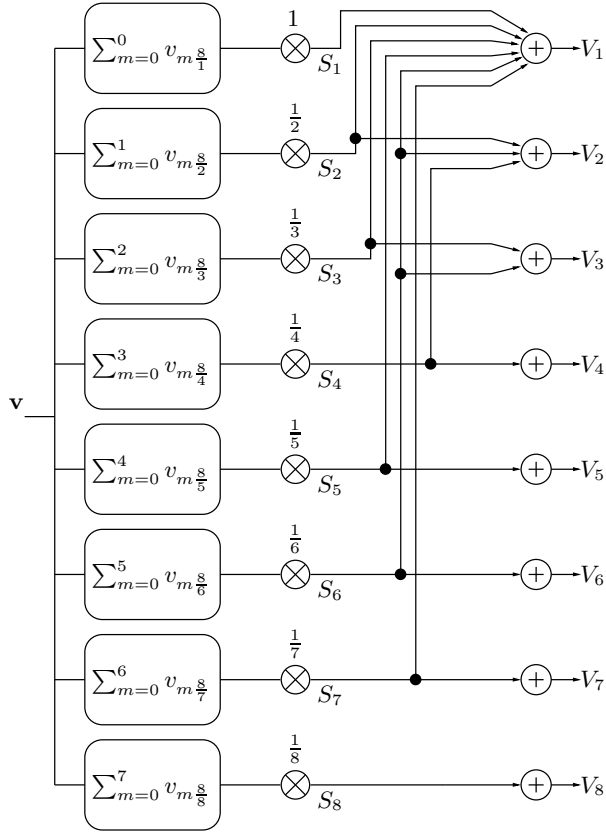}
\caption
{Diagram for computing the AHT for $N=8$.
The boxes compute the averages and the
multipliers implement the scaling operation.
The third layer accounts for the
arithmetic computation based on M\"obius functions.}
\label{diag:tah:8}
\end{figure}

To illustrate the application of above theorem,
consider an 8-point DHT.
Using Theorem~\ref{inv}, Equation~\ref{somaS},
we obtain:
\begin{alignat*}{2}
V_1 & =  S_1-S_2-S_3-S_5+S_6-S_7, & \quad & \\
V_2 & =  S_2-S_4-S_6,             & \quad V_5 & =  S_5, \\
V_3 & =  S_3-S_6,                 & \quad V_6 & =  S_6, \\
V_4 & =  S_4,                     & \quad V_7 & =  S_7.
\end{alignat*}
The component $V_0 = V_8$ can be computed directly from the given samples,
since it represents the mean value of the signal $V_0 = \frac{1}{8} \sum_{m=0}^7v_m$.
Figure~\ref{diag:tah:8} shows a diagram for this computation.

The above theorem and equations completely specified how to compute the discrete Hartley spectrum.
Additionally, the inverse transformation can also be established.
The following result is straightforward.
\begin{corollary}
Inverse discrete Hartley transform components can be computed by
\begin{equation}\label{IAHT}
v_i=
\sum_{l=1}^{\lfloor(N-1)/i \rfloor}
\mu(l)\sigma_{il},
\end{equation}
where
$\sigma_i\triangleq\frac{1}{i} \sum_{m=0}^{i-1}V_{m\frac{N}{i}}$,
$i = 1,\ldots,N-1$.
\endproof
\end{corollary}

The original Arithmetic \emph{Fourier}
Transform had identical
equations to those just derived for the \emph{Hartley} transform
(compare Equation~\ref{sadasiv:9} and Equation~\ref{equacao:48}).
A question arises: since the equations are the same, which spectrum
is actually being evaluated? Fourier or Hartley spectrum?
A clear understanding of underlying arithmetic transform mechanisms will be
possible in the next section.
Once more the reader is asked to put this question
aside for a while, allowing further developments to be derived.

To sum it up, at this point two major questions are accumulated:
(i)~How to handle with fractional indexes? and (ii)~How does same formulae result in different spectra?
Interestingly, both questions had the same answer.

The arithmetic transform algorithm can be summarized in four major steps:
\begin{enumerate}
\item
Index generation, i.e., calculating the indexes of necessary samples (${m\frac{N}{k}}$);

\item
Fractional index samples handling, which requires interpolation;

\item
Computation of averages: $S_k\triangleq\frac{1}{k}\sum_{m=0}^{k-1}v_{m\frac{N}{k}}$;

\item
Computation of spectrum by M{\"o}bius Inversion Formula:
$V_k=\sum_{l=1}^{\lfloor(N-1)/k \rfloor}\mu(l)S_{kl}$.
\end{enumerate}

In the rest of this paper, the step two is addressed.
In the sequel, a mathematical method,
explaining the importance of
the interpolation process
in the arithmetic algorithms,
is derived.

\subsection{Interpolation}

Arithmetic transform theory usually
prescribes
zero- or first-order interpolation
for
the computation of spectral approximations~\cite{Reed:90,Reed:92,Hsu:94}.
In this section, it is shown that
an interpolation process
based on
the known components
(integer index samples)
characterizes
the definition of the
fractional index components, $v_r$, $r\not\in\mathbb{N}$.
This analysis allows a more encompassing
perception of the interpolation mechanisms and gives mathematical tools
for establishing validation constraints to such interpolation process.
In addition,
brief comments on the trade-off
between accuracy and computational cost required by interpolation process
close the section.

\subsubsection{Ideal Interpolation}

What does a fractional index discrete signal component really mean?
The value of $v_r$
for a noninteger value $r$, $r\not\in\mathbb{N}$, can be computed by
\begin{equation}
\label{cas_interp}
\begin{split}
v_r & =  \sum_{k=0}^{N-1}V_{k}\cas\left(\frac{2\pi kr}{N}\right) \\
    & =  \sum_{i=0}^{N-1}v_i\sum_{k=0}^{N-1}
    \cas\left(\frac{2\pi ki}{N}\right)\cas\left(\frac{2\pi kr}{N}\right).
\end{split}
\end{equation}

Defining the  \emph{Hartley weighting function} by
\begin{equation}\label{peso}
w_i(r)\triangleq  \sum_{k=0}^{N-1}
\cas\left(\frac{2\pi ki}{N}\right)
\cas\left(\frac{2\pi kr}{N}\right),
\end{equation}
the value of the signal at fractional indexes can be
found utilizing an $N$-order interpolation expressed by:
\begin{equation}\label{ind_frac}
v_r\triangleq\sum_{i=0}^{N-1}w_i(r)\cdot v_i.
\end{equation}

It is clear that each transform kernel can be associated to
a different weighting function.
Consequently, a different interpolation process for
each weighting function is required.
In the arithmetic transform formalism,
the difference from one transform to another resides in its
interpolation process.

\begin{sloppypar}
The weighting functions satisfy
$\sum_{i=0}^{N-1}{w_i(r)}=1$.
If $r$ is an integer number, then the
orthogonality properties of $\cas(\cdot)$ function~\cite{Bracewell}
make
$w_r(r)=1$ and $w_i(r)=0\ (\forall i\neq r)$.
Therefore, no interpolation is needed.
\end{sloppypar}

\begin{sloppypar}
After some trigonometrical manipulation,
the interpolation weights for several kernels is expressed
by closed formulae.
Let $\Sa(\cdot)$ be the sampling function,
\begin{equation}
\Sa(x)\triangleq
\begin{cases}
{\sin (x)}/{x},& x\neq0, \\
1, & x = 0.
\end{cases}
\end{equation}
\end{sloppypar}

\begin{proposition}\label{closed}
The $N$-point discrete Fourier cosine, Fourier sine, and Hartley transforms
have interpolation weighting functions
given by:

\noindent
\textsl{Fourier cosine Kernel} \\
$
w_i(r)=
\frac{1}{2N}+
\frac{N-1/2}{N}
\left\{
\frac{1}{2}
\frac{\Sa\left(\frac{N-1/2}{N}2\pi(i-r)\right)}{\Sa(\pi(i-r)/N)}
\right.
+
$
\\
\strut\hfill
$
\left.
\frac{1}{2}
\frac{\Sa\left(\frac{N-1/2}{N}2\pi(i+r)\right)}{\Sa(\pi(i+r)/N)}
\right\}.
$

\noindent
\textsl{Fourier sine kernel} \\
$
w_i(r)=
\frac{N-1/2}{N}
\left\{
\frac{1}{2} \frac{\Sa\left(\frac{N-1/2}{N}2\pi(i-r)\right)}{\Sa(\pi(i-r)/N)}
\right.
-
$
\\
\strut\hfill
$
\left.
\frac{1}{2}
\frac{\Sa\left(\frac{N-1/2}{N}2\pi(i+r)\right)}{\Sa(\pi(i+r)/N)}
\right\}.
$

\noindent
\textsl{Hartley kernel} \\
$
w_i(r)=\frac{1}{2N}+
\frac{N-1/2}{N}
\frac{\Sa\left(\frac{N-1/2}{N}2\pi(i-r)\right)}{\Sa(\pi(i-r)/N)}
+
$
\\
\strut\hfill
$
\frac{1}{2N}
\cot \left( \frac{\pi (i+r)}{N}\right)
-
\frac{1}{2N}
\frac{\cos\left(\frac{N-1/2}{N}2\pi(i+r)\right)}{\sin(\pi(i+r)/N)}.
$

\endproof
\end{proposition}

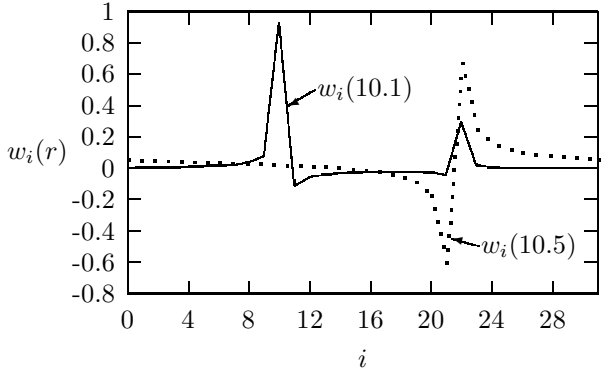
\begin{figure} %
\centering
\setlength{\unitlength}{0.240900pt}
\ifx\plotpoint\undefined\newsavebox{\plotpoint}\fi
\sbox{\plotpoint}{\rule[-0.200pt]{0.400pt}{0.400pt}}%
\begin{picture}(1800,566)(0,0)
\font\gnuplot=cmr10 at 10pt
\gnuplot
\sbox{\plotpoint}{\rule[-0.200pt]{0.400pt}{0.400pt}}%
\put(181.0,123.0){\rule[-0.200pt]{4.818pt}{0.400pt}}
\put(161,123){\makebox(0,0)[r]{-0.8}}
\put(899.0,123.0){\rule[-0.200pt]{4.818pt}{0.400pt}}
\put(181.0,172.0){\rule[-0.200pt]{4.818pt}{0.400pt}}
\put(161,172){\makebox(0,0)[r]{-0.6}}
\put(899.0,172.0){\rule[-0.200pt]{4.818pt}{0.400pt}}
\put(181.0,221.0){\rule[-0.200pt]{4.818pt}{0.400pt}}
\put(161,221){\makebox(0,0)[r]{-0.4}}
\put(899.0,221.0){\rule[-0.200pt]{4.818pt}{0.400pt}}
\put(181.0,271.0){\rule[-0.200pt]{4.818pt}{0.400pt}}
\put(161,271){\makebox(0,0)[r]{-0.2}}
\put(899.0,271.0){\rule[-0.200pt]{4.818pt}{0.400pt}}
\put(181.0,320.0){\rule[-0.200pt]{4.818pt}{0.400pt}}
\put(161,320){\makebox(0,0)[r]{0}}
\put(899.0,320.0){\rule[-0.200pt]{4.818pt}{0.400pt}}
\put(181.0,369.0){\rule[-0.200pt]{4.818pt}{0.400pt}}
\put(161,369){\makebox(0,0)[r]{0.2}}
\put(899.0,369.0){\rule[-0.200pt]{4.818pt}{0.400pt}}
\put(181.0,418.0){\rule[-0.200pt]{4.818pt}{0.400pt}}
\put(161,418){\makebox(0,0)[r]{0.4}}
\put(899.0,418.0){\rule[-0.200pt]{4.818pt}{0.400pt}}
\put(181.0,468.0){\rule[-0.200pt]{4.818pt}{0.400pt}}
\put(161,468){\makebox(0,0)[r]{0.6}}
\put(899.0,468.0){\rule[-0.200pt]{4.818pt}{0.400pt}}
\put(181.0,517.0){\rule[-0.200pt]{4.818pt}{0.400pt}}
\put(161,517){\makebox(0,0)[r]{0.8}}
\put(899.0,517.0){\rule[-0.200pt]{4.818pt}{0.400pt}}
\put(181.0,566.0){\rule[-0.200pt]{4.818pt}{0.400pt}}
\put(161,566){\makebox(0,0)[r]{1}}
\put(899.0,566.0){\rule[-0.200pt]{4.818pt}{0.400pt}}
\put(181.0,123.0){\rule[-0.200pt]{0.400pt}{4.818pt}}
\put(181,82){\makebox(0,0){0}}
\put(181.0,546.0){\rule[-0.200pt]{0.400pt}{4.818pt}}
\put(276.0,123.0){\rule[-0.200pt]{0.400pt}{4.818pt}}
\put(276,82){\makebox(0,0){4}}
\put(276.0,546.0){\rule[-0.200pt]{0.400pt}{4.818pt}}
\put(371.0,123.0){\rule[-0.200pt]{0.400pt}{4.818pt}}
\put(371,82){\makebox(0,0){8}}
\put(371.0,546.0){\rule[-0.200pt]{0.400pt}{4.818pt}}
\put(467.0,123.0){\rule[-0.200pt]{0.400pt}{4.818pt}}
\put(467,82){\makebox(0,0){12}}
\put(467.0,546.0){\rule[-0.200pt]{0.400pt}{4.818pt}}
\put(562.0,123.0){\rule[-0.200pt]{0.400pt}{4.818pt}}
\put(562,82){\makebox(0,0){16}}
\put(562.0,546.0){\rule[-0.200pt]{0.400pt}{4.818pt}}
\put(657.0,123.0){\rule[-0.200pt]{0.400pt}{4.818pt}}
\put(657,82){\makebox(0,0){20}}
\put(657.0,546.0){\rule[-0.200pt]{0.400pt}{4.818pt}}
\put(752.0,123.0){\rule[-0.200pt]{0.400pt}{4.818pt}}
\put(752,82){\makebox(0,0){24}}
\put(752.0,546.0){\rule[-0.200pt]{0.400pt}{4.818pt}}
\put(848.0,123.0){\rule[-0.200pt]{0.400pt}{4.818pt}}
\put(848,82){\makebox(0,0){28}}
\put(848.0,546.0){\rule[-0.200pt]{0.400pt}{4.818pt}}
\put(181.0,123.0){\rule[-0.200pt]{177.784pt}{0.400pt}}
\put(919.0,123.0){\rule[-0.200pt]{0.400pt}{106.719pt}}
\put(181.0,566.0){\rule[-0.200pt]{177.784pt}{0.400pt}}
\put(40,344){\makebox(0,0){$w_i(r)$}}
\put(550,21){\makebox(0,0){$i$}}
\put(479,443){\makebox(0,0)[l]{$w_{i}(10.1)$}}
\put(736,197){\makebox(0,0)[l]{$w_{i}(10.5)$}}
\put(181.0,123.0){\rule[-0.200pt]{0.400pt}{106.719pt}}
\multiput(472.73,441.92)(-0.864,-0.497){47}{\rule{0.788pt}{0.120pt}}
\multiput(474.36,442.17)(-41.364,-25.000){2}{\rule{0.394pt}{0.400pt}}
\put(433,418){\vector(-2,-1){0}}
\multiput(726.63,197.58)(-1.832,0.492){21}{\rule{1.533pt}{0.119pt}}
\multiput(729.82,196.17)(-39.817,12.000){2}{\rule{0.767pt}{0.400pt}}
\put(690,209){\vector(-4,1){0}}
\put(181,320){\usebox{\plotpoint}}
\put(181,319.67){\rule{5.782pt}{0.400pt}}
\multiput(181.00,319.17)(12.000,1.000){2}{\rule{2.891pt}{0.400pt}}
\put(252,320.67){\rule{5.782pt}{0.400pt}}
\multiput(252.00,320.17)(12.000,1.000){2}{\rule{2.891pt}{0.400pt}}
\put(276,321.67){\rule{5.782pt}{0.400pt}}
\multiput(276.00,321.17)(12.000,1.000){2}{\rule{2.891pt}{0.400pt}}
\put(300,322.67){\rule{5.782pt}{0.400pt}}
\multiput(300.00,322.17)(12.000,1.000){2}{\rule{2.891pt}{0.400pt}}
\put(324,323.67){\rule{5.782pt}{0.400pt}}
\multiput(324.00,323.17)(12.000,1.000){2}{\rule{2.891pt}{0.400pt}}
\multiput(348.00,325.60)(3.259,0.468){5}{\rule{2.400pt}{0.113pt}}
\multiput(348.00,324.17)(18.019,4.000){2}{\rule{1.200pt}{0.400pt}}
\multiput(371.00,329.58)(1.225,0.491){17}{\rule{1.060pt}{0.118pt}}
\multiput(371.00,328.17)(21.800,10.000){2}{\rule{0.530pt}{0.400pt}}
\multiput(395.58,339.00)(0.496,4.412){45}{\rule{0.120pt}{3.583pt}}
\multiput(394.17,339.00)(24.000,201.563){2}{\rule{0.400pt}{1.792pt}}
\multiput(419.58,529.87)(0.496,-5.406){45}{\rule{0.120pt}{4.367pt}}
\multiput(418.17,538.94)(24.000,-246.937){2}{\rule{0.400pt}{2.183pt}}
\multiput(443.00,292.58)(0.864,0.494){25}{\rule{0.786pt}{0.119pt}}
\multiput(443.00,291.17)(22.369,14.000){2}{\rule{0.393pt}{0.400pt}}
\multiput(467.00,306.60)(3.259,0.468){5}{\rule{2.400pt}{0.113pt}}
\multiput(467.00,305.17)(18.019,4.000){2}{\rule{1.200pt}{0.400pt}}
\put(490,310.17){\rule{4.900pt}{0.400pt}}
\multiput(490.00,309.17)(13.830,2.000){2}{\rule{2.450pt}{0.400pt}}
\put(514,311.67){\rule{5.782pt}{0.400pt}}
\multiput(514.00,311.17)(12.000,1.000){2}{\rule{2.891pt}{0.400pt}}
\put(538,312.67){\rule{5.782pt}{0.400pt}}
\multiput(538.00,312.17)(12.000,1.000){2}{\rule{2.891pt}{0.400pt}}
\put(205.0,321.0){\rule[-0.200pt]{11.322pt}{0.400pt}}
\put(633,312.67){\rule{5.782pt}{0.400pt}}
\multiput(633.00,313.17)(12.000,-1.000){2}{\rule{2.891pt}{0.400pt}}
\multiput(657.00,311.94)(3.406,-0.468){5}{\rule{2.500pt}{0.113pt}}
\multiput(657.00,312.17)(18.811,-4.000){2}{\rule{1.250pt}{0.400pt}}
\multiput(681.58,309.00)(0.496,1.746){45}{\rule{0.120pt}{1.483pt}}
\multiput(680.17,309.00)(24.000,79.921){2}{\rule{0.400pt}{0.742pt}}
\multiput(705.58,386.88)(0.496,-1.429){45}{\rule{0.120pt}{1.233pt}}
\multiput(704.17,389.44)(24.000,-65.440){2}{\rule{0.400pt}{0.617pt}}
\multiput(729.00,322.95)(4.927,-0.447){3}{\rule{3.167pt}{0.108pt}}
\multiput(729.00,323.17)(16.427,-3.000){2}{\rule{1.583pt}{0.400pt}}
\put(752,319.67){\rule{5.782pt}{0.400pt}}
\multiput(752.00,320.17)(12.000,-1.000){2}{\rule{2.891pt}{0.400pt}}
\put(562.0,314.0){\rule[-0.200pt]{17.104pt}{0.400pt}}
\put(776.0,320.0){\rule[-0.200pt]{34.449pt}{0.400pt}}
\sbox{\plotpoint}{\rule[-0.500pt]{1.000pt}{1.000pt}}%
\put(181,332){\usebox{\plotpoint}}
\multiput(181,332)(20.738,-0.864){2}{\usebox{\plotpoint}}
\put(222.48,330.27){\usebox{\plotpoint}}
\put(243.22,330.00){\usebox{\plotpoint}}
\put(263.97,329.50){\usebox{\plotpoint}}
\put(284.71,328.64){\usebox{\plotpoint}}
\put(305.44,327.77){\usebox{\plotpoint}}
\multiput(324,327)(20.738,-0.864){2}{\usebox{\plotpoint}}
\put(367.67,326.00){\usebox{\plotpoint}}
\put(388.41,325.27){\usebox{\plotpoint}}
\put(409.15,324.41){\usebox{\plotpoint}}
\put(429.89,323.55){\usebox{\plotpoint}}
\put(450.63,322.68){\usebox{\plotpoint}}
\put(471.36,321.81){\usebox{\plotpoint}}
\multiput(490,321)(20.684,-1.724){2}{\usebox{\plotpoint}}
\put(533.46,317.38){\usebox{\plotpoint}}
\put(554.15,315.65){\usebox{\plotpoint}}
\put(574.70,312.88){\usebox{\plotpoint}}
\put(595.10,309.10){\usebox{\plotpoint}}
\put(615.16,303.98){\usebox{\plotpoint}}
\multiput(633,297)(15.620,-13.668){2}{\usebox{\plotpoint}}
\multiput(657,276)(4.625,-20.234){5}{\usebox{\plotpoint}}
\multiput(681,171)(1.587,20.695){15}{\usebox{\plotpoint}}
\multiput(705,484)(4.625,-20.234){5}{\usebox{\plotpoint}}
\multiput(729,379)(15.328,-13.995){2}{\usebox{\plotpoint}}
\put(762.42,354.09){\usebox{\plotpoint}}
\put(782.12,347.72){\usebox{\plotpoint}}
\multiput(800,344)(20.473,-3.412){2}{\usebox{\plotpoint}}
\put(843.61,338.37){\usebox{\plotpoint}}
\put(864.29,336.58){\usebox{\plotpoint}}
\put(885.00,335.42){\usebox{\plotpoint}}
\put(905.71,334.11){\usebox{\plotpoint}}
\put(919,333){\usebox{\plotpoint}}
\end{picture}
\vspace{-.8cm}
\caption{Hartley weighting functions used to interpolate $v_{10.1}$ and $v_{10.5}$ ($N=32$ blocklength).}
\label{fig:pesos}
\end{figure}

\begin{figure}[!t]
\centering
\subfigure[Fourier cosine kernel]{\epsfig{file=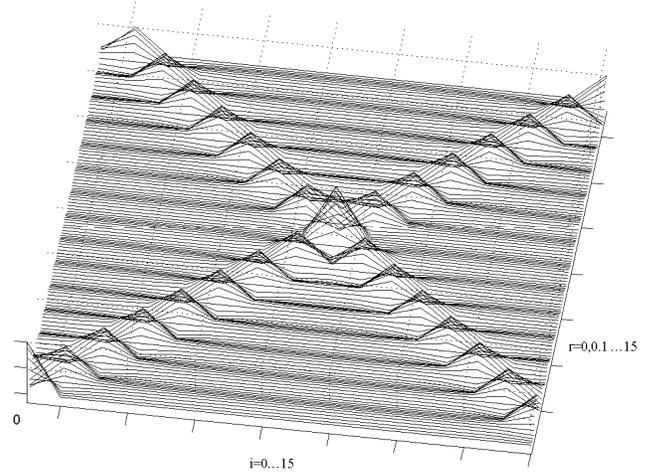,width=.95\linewidth}}\\
\subfigure[Hartley kernel]{\epsfig{file=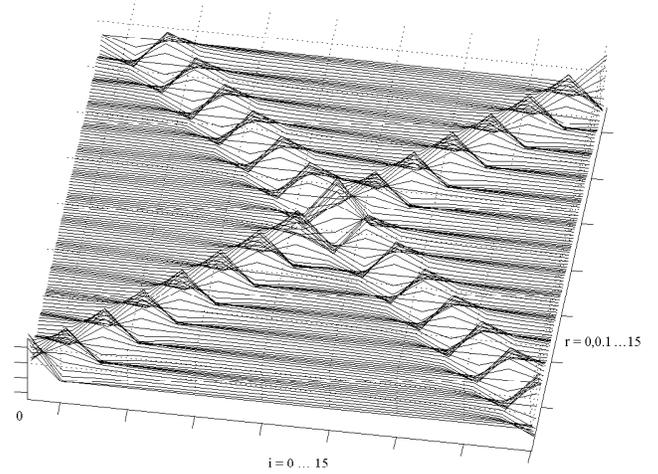,width=.95\linewidth}}
\caption{
These curve families represent the
weighting profile~$w_i(r)$ for the
Fourier cosine kernel~(a) and
Hartley kernel~(b).
The parameters
$N=16$, $r = 0,0.1\ldots15$, and $i = 0\ldots15$ are employed.
Observe that $i\in \mathbb{N}$ and $r \in\mathbb{R}$.
The maximum values are achieved at $i =0$ or $i = N/2 = 8$ (central peak)
and they corresponds to the unity.
The value of the local maxima are exactly $1/2$, occurring when $r$ is integer.
Note that each curve is almost null everywhere, except at the vicinities of $i\approx r$ or $i \approx N-r$.
}
\label{fig:profile}
\end{figure}

To exemplify, Figure~\ref{fig:pesos} shows two weighting functions used
to compute $v_{10.1}$ and $v_{10.5}$ in the arithmetic Hartley transform.
Figure~\ref{fig:profile} shows the weighting profile for the Fourier cosine and Hartley kernels.
These functions are calculated by closed formulae.

With this proposition, the
mathematical description of the AHT algorithm
is completed.
The derived formulae furnish
the \emph{exact} value of the spectral components.
On the other hand, the computational complexity of the ideal interpolation
implementation is similar to
the direct implementation, i.e., computing the transform by its plain definition:
$V_k=\frac{1}{N}\sum_{i=0}^{N-1}v_i\cas\left(\frac{2\pi}{N}ki\right)$.
To address this issue, a non-ideal interpolation scheme is proposed.

\subsubsection{Non-Ideal Interpolation}

According to the index generation (${m\frac{N}{k}}$),
the number $R$ of points that require interpolation is upper bounded by
$
R\leq\sum_{d \nmid N}{d-1}
$.
This sum represents the number of samples with fractional index.
Consequently, this approach can be attractive for large non-prime blocklength~$N$
with great number of factors, because it would require a smaller
number of interpolations.

The next task is to find simpler formulae for the weighting functions,
assuming large blocklength.
Instead of using the exact weighting functions,
the limit $N\!\!\to\!\infty$ is examined and
used to derive asymptotic approximations of the weighting function.

\begin{proposition}\label{Asymptotic}
A continuous approximation for the interpolation weighting function for sufficiently large $N$ is given by:

\noindent
\textsl{Fourier cosine kernel} \\
\strut\hfill
$
\hat{w}_i(r)\approx\frac{ \Sa(2\pi(i-r))}{2}+\frac{\Sa(2\pi(i+r)) }{2}
$
\hfill\strut

\noindent
\textsl{Fourier sine kernel} \\
\strut\hfill
$
\hat{w}_i(r)\approx \frac{\Sa(2\pi(i-r))}{2}-\frac{\Sa(2\pi(i+r)) }{2}
$
\hfill\strut

\noindent
\textsl{Hartley kernel} \\
\strut\hfill
$
\hat{w}_i(r)\approx
\Sa(2\pi(i-r))+
\frac{1-\cos 2\pi r }{2\pi(i+r)}.
$
\hfill\strut
\end{proposition}
\endproof

In terms of the Hilbert transform and
the $\Sa(\cdot)$ function,
the asymptotic weighting function for Hartley kernel is given by
\begin{align}
\hat{w}_i(r)\approx &
\Sa(2\pi(i-r))-
\mathcal{H}il
\Big\{
\Sa(2\pi(i+r))
\Big\},\\
\intertext{or alternatively,}
\hat{w}_i(r)\approx &
\Sa(2\pi(i-r))
-
\Ca(2\pi(i+r))
-\\ \nonumber
&
\mathcal{H}il
\Big\{
\delta\big(2(r+i)\big)
\Big\},
\end{align}
where ${\cal H}il$ denotes the Hilbert transform,
$\Ca(x)\triangleq\frac{\cos x}{x}$, for $x\neq0$, is the co-sampling
function and $\delta(x)$ is the Dirac impulse.

\paragraph{Zero-order Interpolation.}
Rounding the fractional index provides
the zero-order interpolation.
The estimated (interpolated) signal $\hat{v}_j$
is then expressed by
$\hat{v}_j=v_{[j]}$, where~$[\cdot]$ is a function
which rounds off its argument to its nearest integer.
Examining the asymptotic behavior of the weighting function
for Fourier cosine kernel,
we derive the following results:
\begin{equation}\label{zero=cos}
\begin{split}
\hat{w}_i(r)&\approx 0,  \quad \forall\; i\neq {[r],N-[r]}, \\
\hat{w}_{[r]}(r)& \approx \frac{\Sa\left( 2\pi([r]-r)\right)}{2}  \approx
\frac{1}{2}, \\
\hat{w}_{N-[r]}(r)& \approx \frac{\Sa \left( 2\pi([r]-r)\right)}{2}
\approx\frac{1}{2}.\\
\end{split}
\end{equation}
Under the above assumptions, the Equation~\ref{ind_frac} furnishes:
\begin{equation}
\begin{split}
\hat{v}_r &\approx w_{[r]}(r)v_{[r]}+w_{N-[r]}(r)v_{N-[r]}\\
& \approx \frac{1}{2}v_{[r]} + \frac{1}{2}v_{N-[r]}.\\
\end{split}
\end{equation}

Thus, for even signal ($v_k=v_{N-k}$) ,
the approximated
value
of the interpolated sample
is roughly given by $\hat{v}_r \approx v_{[r]}$.
It is straightforward to see that the influence of odd part of the signal
vanishes in the zero-order interpolation.
Besides zero-order interpolation is ``blind'' to the odd component of a signal.
In fact, as show by the set of Equations~\ref{zero=cos},
zero-order interpolation
is an (indeed good) approximation to the cosine asymptotic weighting function.
This puts some light on the role of the zero-order interpolations and
its relation to the earliest versions of the arithmetic Fourier transform,
which can not analyze odd periodic signals.

Zero-order interpolation
has been intuitively employed in previous work by Tufts, Reed {\em et alli}~\cite{Reed:90,Tufts:Sadasiv:88a,Reed:92}.
Hsu, in his Ph.D. dissertation, derives an analysis of first-order
interpolation effect~\cite{Hsu:94}.

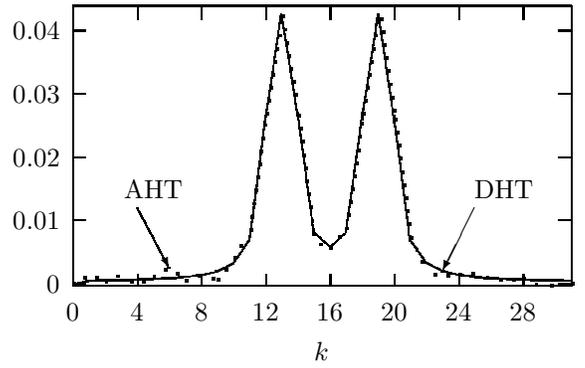
\begin{figure} %
\centering
\setlength{\unitlength}{0.240900pt}
\ifx\plotpoint\undefined\newsavebox{\plotpoint}\fi
\sbox{\plotpoint}{\rule[-0.200pt]{0.400pt}{0.400pt}}%
\begin{picture}(944,603)(0,0)
\font\gnuplot=cmr10 at 10pt
\gnuplot
\sbox{\plotpoint}{\rule[-0.200pt]{0.400pt}{0.400pt}}%
\put(140.0,126.0){\rule[-0.200pt]{4.818pt}{0.400pt}}
\put(120,126){\makebox(0,0)[r]{0}}
\put(903.0,126.0){\rule[-0.200pt]{4.818pt}{0.400pt}}
\put(140.0,226.0){\rule[-0.200pt]{4.818pt}{0.400pt}}
\put(120,226){\makebox(0,0)[r]{0.01}}
\put(903.0,226.0){\rule[-0.200pt]{4.818pt}{0.400pt}}
\put(140.0,325.0){\rule[-0.200pt]{4.818pt}{0.400pt}}
\put(120,325){\makebox(0,0)[r]{0.02}}
\put(903.0,325.0){\rule[-0.200pt]{4.818pt}{0.400pt}}
\put(140.0,425.0){\rule[-0.200pt]{4.818pt}{0.400pt}}
\put(120,425){\makebox(0,0)[r]{0.03}}
\put(903.0,425.0){\rule[-0.200pt]{4.818pt}{0.400pt}}
\put(140.0,524.0){\rule[-0.200pt]{4.818pt}{0.400pt}}
\put(120,524){\makebox(0,0)[r]{0.04}}
\put(903.0,524.0){\rule[-0.200pt]{4.818pt}{0.400pt}}
\put(140.0,123.0){\rule[-0.200pt]{0.400pt}{4.818pt}}
\put(140,82){\makebox(0,0){0}}
\put(140.0,543.0){\rule[-0.200pt]{0.400pt}{4.818pt}}
\put(241.0,123.0){\rule[-0.200pt]{0.400pt}{4.818pt}}
\put(241,82){\makebox(0,0){4}}
\put(241.0,543.0){\rule[-0.200pt]{0.400pt}{4.818pt}}
\put(342.0,123.0){\rule[-0.200pt]{0.400pt}{4.818pt}}
\put(342,82){\makebox(0,0){8}}
\put(342.0,543.0){\rule[-0.200pt]{0.400pt}{4.818pt}}
\put(443.0,123.0){\rule[-0.200pt]{0.400pt}{4.818pt}}
\put(443,82){\makebox(0,0){12}}
\put(443.0,543.0){\rule[-0.200pt]{0.400pt}{4.818pt}}
\put(544.0,123.0){\rule[-0.200pt]{0.400pt}{4.818pt}}
\put(544,82){\makebox(0,0){16}}
\put(544.0,543.0){\rule[-0.200pt]{0.400pt}{4.818pt}}
\put(645.0,123.0){\rule[-0.200pt]{0.400pt}{4.818pt}}
\put(645,82){\makebox(0,0){20}}
\put(645.0,543.0){\rule[-0.200pt]{0.400pt}{4.818pt}}
\put(746.0,123.0){\rule[-0.200pt]{0.400pt}{4.818pt}}
\put(746,82){\makebox(0,0){24}}
\put(746.0,543.0){\rule[-0.200pt]{0.400pt}{4.818pt}}
\put(847.0,123.0){\rule[-0.200pt]{0.400pt}{4.818pt}}
\put(847,82){\makebox(0,0){28}}
\put(847.0,543.0){\rule[-0.200pt]{0.400pt}{4.818pt}}
\put(140.0,123.0){\rule[-0.200pt]{188.625pt}{0.400pt}}
\put(923.0,123.0){\rule[-0.200pt]{0.400pt}{105.996pt}}
\put(140.0,563.0){\rule[-0.200pt]{188.625pt}{0.400pt}}
\put(531,21){\makebox(0,0){$k$}}
\put(771,275){\makebox(0,0)[l]{DHT}}
\put(221,275){\makebox(0,0)[l]{AHT}}
\put(140.0,123.0){\rule[-0.200pt]{0.400pt}{105.996pt}}
\put(771,246){\vector(-1,-2){50}}
\multiput(241.58,242.49)(0.498,-0.934){99}{\rule{0.120pt}{0.845pt}}
\multiput(240.17,244.25)(51.000,-93.246){2}{\rule{0.400pt}{0.423pt}}
\put(292,151){\vector(1,-2){0}}
\put(140,126){\usebox{\plotpoint}}
\multiput(140.00,126.59)(2.714,0.477){7}{\rule{2.100pt}{0.115pt}}
\multiput(140.00,125.17)(20.641,5.000){2}{\rule{1.050pt}{0.400pt}}
\put(165,130.67){\rule{6.263pt}{0.400pt}}
\multiput(165.00,130.17)(13.000,1.000){2}{\rule{3.132pt}{0.400pt}}
\put(216,131.67){\rule{6.023pt}{0.400pt}}
\multiput(216.00,131.17)(12.500,1.000){2}{\rule{3.011pt}{0.400pt}}
\put(241,132.67){\rule{6.023pt}{0.400pt}}
\multiput(241.00,132.17)(12.500,1.000){2}{\rule{3.011pt}{0.400pt}}
\put(266,133.67){\rule{6.263pt}{0.400pt}}
\multiput(266.00,133.17)(13.000,1.000){2}{\rule{3.132pt}{0.400pt}}
\put(292,135.17){\rule{5.100pt}{0.400pt}}
\multiput(292.00,134.17)(14.415,2.000){2}{\rule{2.550pt}{0.400pt}}
\multiput(317.00,137.61)(5.374,0.447){3}{\rule{3.433pt}{0.108pt}}
\multiput(317.00,136.17)(17.874,3.000){2}{\rule{1.717pt}{0.400pt}}
\multiput(342.00,140.59)(2.208,0.482){9}{\rule{1.767pt}{0.116pt}}
\multiput(342.00,139.17)(21.333,6.000){2}{\rule{0.883pt}{0.400pt}}
\multiput(367.00,146.58)(1.012,0.493){23}{\rule{0.900pt}{0.119pt}}
\multiput(367.00,145.17)(24.132,13.000){2}{\rule{0.450pt}{0.400pt}}
\multiput(393.58,159.00)(0.497,0.722){47}{\rule{0.120pt}{0.676pt}}
\multiput(392.17,159.00)(25.000,34.597){2}{\rule{0.400pt}{0.338pt}}
\multiput(418.58,195.00)(0.497,3.847){47}{\rule{0.120pt}{3.140pt}}
\multiput(417.17,195.00)(25.000,183.483){2}{\rule{0.400pt}{1.570pt}}
\multiput(443.58,385.00)(0.497,3.340){47}{\rule{0.120pt}{2.740pt}}
\multiput(442.17,385.00)(25.000,159.313){2}{\rule{0.400pt}{1.370pt}}
\multiput(468.58,539.37)(0.497,-3.112){49}{\rule{0.120pt}{2.562pt}}
\multiput(467.17,544.68)(26.000,-154.683){2}{\rule{0.400pt}{1.281pt}}
\multiput(494.58,377.43)(0.497,-3.705){47}{\rule{0.120pt}{3.028pt}}
\multiput(493.17,383.72)(25.000,-176.715){2}{\rule{0.400pt}{1.514pt}}
\multiput(519.00,205.92)(0.542,-0.496){43}{\rule{0.535pt}{0.120pt}}
\multiput(519.00,206.17)(23.890,-23.000){2}{\rule{0.267pt}{0.400pt}}
\multiput(544.00,184.58)(0.542,0.496){43}{\rule{0.535pt}{0.120pt}}
\multiput(544.00,183.17)(23.890,23.000){2}{\rule{0.267pt}{0.400pt}}
\multiput(569.58,207.00)(0.497,3.560){49}{\rule{0.120pt}{2.915pt}}
\multiput(568.17,207.00)(26.000,176.949){2}{\rule{0.400pt}{1.458pt}}
\multiput(595.58,390.00)(0.497,3.238){47}{\rule{0.120pt}{2.660pt}}
\multiput(594.17,390.00)(25.000,154.479){2}{\rule{0.400pt}{1.330pt}}
\multiput(620.58,538.63)(0.497,-3.340){47}{\rule{0.120pt}{2.740pt}}
\multiput(619.17,544.31)(25.000,-159.313){2}{\rule{0.400pt}{1.370pt}}
\multiput(645.58,371.97)(0.497,-3.847){47}{\rule{0.120pt}{3.140pt}}
\multiput(644.17,378.48)(25.000,-183.483){2}{\rule{0.400pt}{1.570pt}}
\multiput(670.58,192.29)(0.497,-0.693){49}{\rule{0.120pt}{0.654pt}}
\multiput(669.17,193.64)(26.000,-34.643){2}{\rule{0.400pt}{0.327pt}}
\multiput(696.00,157.92)(0.972,-0.493){23}{\rule{0.869pt}{0.119pt}}
\multiput(696.00,158.17)(23.196,-13.000){2}{\rule{0.435pt}{0.400pt}}
\multiput(721.00,144.93)(2.208,-0.482){9}{\rule{1.767pt}{0.116pt}}
\multiput(721.00,145.17)(21.333,-6.000){2}{\rule{0.883pt}{0.400pt}}
\multiput(746.00,138.95)(5.374,-0.447){3}{\rule{3.433pt}{0.108pt}}
\multiput(746.00,139.17)(17.874,-3.000){2}{\rule{1.717pt}{0.400pt}}
\put(771,135.17){\rule{5.300pt}{0.400pt}}
\multiput(771.00,136.17)(15.000,-2.000){2}{\rule{2.650pt}{0.400pt}}
\put(797,133.67){\rule{6.023pt}{0.400pt}}
\multiput(797.00,134.17)(12.500,-1.000){2}{\rule{3.011pt}{0.400pt}}
\put(822,132.67){\rule{6.023pt}{0.400pt}}
\multiput(822.00,133.17)(12.500,-1.000){2}{\rule{3.011pt}{0.400pt}}
\put(847,131.67){\rule{6.023pt}{0.400pt}}
\multiput(847.00,132.17)(12.500,-1.000){2}{\rule{3.011pt}{0.400pt}}
\put(191.0,132.0){\rule[-0.200pt]{6.022pt}{0.400pt}}
\put(898,130.67){\rule{6.023pt}{0.400pt}}
\multiput(898.00,131.17)(12.500,-1.000){2}{\rule{3.011pt}{0.400pt}}
\put(872.0,132.0){\rule[-0.200pt]{6.263pt}{0.400pt}}
\sbox{\plotpoint}{\rule[-0.500pt]{1.000pt}{1.000pt}}%
\put(140,126){\usebox{\plotpoint}}
\multiput(140,126)(18.109,10.141){2}{\usebox{\plotpoint}}
\put(176.84,134.99){\usebox{\plotpoint}}
\multiput(191,129)(18.712,8.982){2}{\usebox{\plotpoint}}
\put(231.85,130.22){\usebox{\plotpoint}}
\put(249.45,128.73){\usebox{\plotpoint}}
\multiput(266,138)(18.564,9.282){2}{\usebox{\plotpoint}}
\put(303.11,142.11){\usebox{\plotpoint}}
\multiput(317,131)(19.529,7.030){2}{\usebox{\plotpoint}}
\put(359.06,134.54){\usebox{\plotpoint}}
\multiput(367,132)(12.377,16.661){2}{\usebox{\plotpoint}}
\multiput(393,167)(10.253,18.046){2}{\usebox{\plotpoint}}
\multiput(418,211)(2.825,20.562){9}{\usebox{\plotpoint}}
\multiput(443,393)(3.347,20.484){8}{\usebox{\plotpoint}}
\multiput(468,546)(3.412,-20.473){7}{\usebox{\plotpoint}}
\multiput(494,390)(2.587,-20.594){10}{\usebox{\plotpoint}}
\put(528.01,188.12){\usebox{\plotpoint}}
\multiput(544,183)(13.287,15.945){2}{\usebox{\plotpoint}}
\multiput(569,213)(2.967,20.542){9}{\usebox{\plotpoint}}
\multiput(595,393)(3.020,20.535){8}{\usebox{\plotpoint}}
\multiput(620,563)(2.871,-20.556){9}{\usebox{\plotpoint}}
\multiput(645,384)(2.765,-20.571){9}{\usebox{\plotpoint}}
\multiput(670,198)(8.870,-18.764){3}{\usebox{\plotpoint}}
\put(707.96,140.61){\usebox{\plotpoint}}
\put(728.43,138.59){\usebox{\plotpoint}}
\multiput(746,140)(20.689,1.655){2}{\usebox{\plotpoint}}
\put(789.70,136.25){\usebox{\plotpoint}}
\put(810.02,132.44){\usebox{\plotpoint}}
\put(830.63,132.04){\usebox{\plotpoint}}
\multiput(847,134)(18.998,-8.359){2}{\usebox{\plotpoint}}
\put(890.47,123.00){\usebox{\plotpoint}}
\put(911.13,124.58){\usebox{\plotpoint}}
\put(923,126){\usebox{\plotpoint}}
\end{picture}
\vspace{-.3cm}
\caption{The discrete Hartley transform computed
by definition (solid line) and by arithmetic transform algorithm (dotted line).
Simulation data: $f(t) = \cos(90\pi t)\left(  t - \frac{1}{2} \right)^2$, $t=0\ldots1$, $N=32$.}
\label{aht:dht}
\end{figure}

\paragraph{Interpolation Order.}

Let $M_m$ be a set with the $m$ ($<N$) most significant coefficients $w_i(r)$.
For zero-order interpolation, $m=1$.
Increasing the value of $m$, the interpolation process would gradually be
improved, because more coefficients would be retained.
Proceeding in this way,
the following calculation
performs
a non-ideal interpolation:
\begin{equation}
\hat{v}_r=
\frac{1}{\eta}
\sum_{i \in M_m}
w_i(r)
\cdot v_i,
\end{equation}
where $\eta \triangleq \sum_{j \in M_m}w_j(r)$ is a normalization factor.

Figure~\ref{aht:dht} presents a 32-point discrete Hartley transform of
a particular signal computed
by the plain definition of the DHT and by the arithmetic transform method using $m=2$.

\section{Conclusions}

This paper supplied a short survey on arithmetic transforms.
The arithmetic Fourier transform is explained and its three main 1-D
versions were described.
Furthermore, some comments on implementation,
challenging points, and advantages
of the arithmetic transforms were discussed via simple examples.

In this paper, the introduction of the AHT emphasized the
key point of the arithmetic transforms:
the interpolation process.
It was shown that the fundamental equations of the arithmetic transform algorithms
were essentially the same (regardless the kernel).
This property could open path to
the implementation of ``universal transformers''.
In this type of construction, the circuitry would be the same
for several transforms,
except the interpolation module.
A different interpolation module would reflect different
transform (Fourier, Hartley, Fourier cosine).
Finally, this paper could be taken as a starting point to those who want to
investigate arithmetic transforms.

\section{Acknowledgments}

The first author would like to thank
Emeritus Professor Dr.~Irving~S. Reed,
University of Southern California,
for kindly sending him a copy of
Chin-Chi~Hsu's Ph.D. thesis~\cite{Hsu:94}.

{\small
\bibliographystyle{IEEEtran}
\bibliography{bib_aft_clean}
}

\strut

\vspace{1cm}

\appendix

\section{Additional References}

The following references are
not primarily concerned with arithmetic transforms,
however they are relevant.

\def\refname{}
\strut

\end{document}